\documentclass[openany, a4paper, 12pt]{article}
\usepackage{amsfonts}
\usepackage{fancyhdr}
\usepackage{amssymb}
 \usepackage{mathrsfs,amsfonts,amsmath}
 \usepackage{color}

\makeatletter

\@addtoreset{equation}{section}
\makeatother
\newtheorem{Theorem}{Theorem}[section]
\newtheorem{Remark}[Theorem]{Remark}
\newtheorem{Lemma}[Theorem]{Lemma}

\begin{document}
\title{\textbf{Polynomial and exponential stability of $\theta$-EM approximations to
a class of stochastic differential equations}}
\author{Yunjiao Hu$^a$,\ Guangqiang Lan$^a$\footnote{Corresponding author. Email:
langq@mail.buct.edu.cn. Supported by
China Scholarship Council, National Natural Science Foundation of
China (NSFC11026142) and Beijing Higher Education Young Elite Teacher Project (YETP0516).},\  Chong Zhang$^a$
\\ \small $^{a}$School of Science, Beijing University of Chemical Technology, Beijing 100029, China}

\date{}

\maketitle

\begin{abstract}
Both the mean square polynomial stability and exponential stability
of $\theta$ Euler-Maruyama approximation solutions of stochastic
differential equations will be investigated for each $0\le\theta\le
1$ by using an auxiliary function $F$ (see the following definition
(\ref{dingyi})). Sufficient conditions are obtained to ensure the
polynomial and exponential stability of the numerical
approximations. The results in Liu et al \cite{LFM} will be improved
and generalized to more general cases. Several examples and non
stability results are presented to support our conclusions.
\end{abstract}



\section{Introduction}

\noindent

Given a probability space $(\Omega,\mathscr{F},P)$ endowed with a
complete filtration $(\mathscr{F}_t)_{t\geq 0}$. Let
$d,m\in\mathbb{N}$ be arbitrarily fixed. We consider the following
stochastic differential equations (SDEs)
\begin{equation}\label{sde}dX_t=f(X_t,t)dt+g(X_t,t)dB_t,\
X_0=x_0\in \mathbb{R}^d,  \end{equation}

\noindent where the initial $x_0\in \mathbb{R}^d, (B_t)_{t\geq0}$ is
an $m$-dimensional standard $\mathscr{F}_t$-Brownian motion,
$f:(t,x)\in[0,\infty)\times\mathbb{R}^d\mapsto f(t,x)\in\mathbb{R}^d$
and $g:(t,x)\in[0,\infty)\times\mathbb{R}^d\mapsto\sigma(t,x)\in
\mathbb{R}^d\otimes\mathbb{R}^m$ are both Borel measurable functions.

The corresponding $\theta$ Euler-Maruyama ($\theta$-EM) approximation
(or the so called stochastic theta method) of the above SDE is
\begin{equation}\label{SEM} X_{k+1}=X_k+[(1-\theta)f(X_k,k\Delta t)+
\theta f(X_{k+1},(k+1)\Delta t)]\Delta t+g(X_k,k\Delta t)\Delta B_k,\end{equation}

\noindent where $X_0:=x_0$, $\Delta t$ is a constant step size,
$\theta\in [0,1]$ is a fixed parameter, $\Delta B_k:=B((k+1)\Delta
t)-B(k\Delta t)$ is the increment of Brownian motion. Note that
$\theta$-EM includes the classical EM method ($\theta=0$), the
backward EM method ($\theta=1$) and the so-called trapezoidal method
($\theta=\frac{1}{2}$).

Throughout of this paper, we simply assume that the coefficients $f$
and $g$ satisfy the following local Lipschitz condition:

For every integer $r\ge1$ and any $t\ge0$, there exists a positive
constant $\bar{K}_{r,t}$ such that for any $x,y\in\mathbb{R}^d$ with
$\max\{|x|,|y|\}\le r,$
\begin{equation}\label{local}
\max\{|f(x,t)-f(y,t)|,|g(x,t)-g(y,t)|\}\le \bar{K}_{r,t}|x-y|.
\end{equation}

Condition (\ref{local}) could make sure that equation (\ref{sde})
has a unique solution, which is denoted by
$X_t(x_0)\in\mathbb{R}^d,$ (this condition could be weakened to more
generalized condition, see e.g. \cite{Lan,Lan1}).

Stability theory is one of the central problems in numerical
analysis. The stability concepts of numerical approximation for SDEs
mainly include moment stability (M-stability) and almost sure
stability (trajectory stability). Results concerned with different
kinds of stability analysis for numerical methods can be found in
many literatures.

For example, Baker and Buckwar \cite{BB} dealt with the $p$-th
moment exponential stability of stochastic delay differential
equations when the coefficients are both globally Lipschitz
continuous, Higham \cite{Higham1,Higham2} considered the scalar
linear case and Higham et al. \cite{HMS} for one sided Lipschitz and
the linear growth condition. Other results concerned with moment
stability can be found in the Mao's monograph \cite{Mao}, Higham et
al \cite{HMY}, Zong et al \cite{ZW}, Pang et al \cite{PDM}, Szpruch
\cite{Szpruch} (for the so called $V$-stability) and references
therein.

For the almost sure stability of numerical approximation for SDEs,
by Borel-Cantelli lemma and Chebyshev inequality, recently, Wu et al
\cite{WMS} investgated the almost sure exponential stability of the
stochastic theta method by the continuous and discrete semi
martingale convergence theorems (see Rodkina and Schurz \cite{RS}
for details), Chen and Wu \cite{CW} and Mao and Szpruch \cite{MS}
also used the same method to prove the almost sure stability of the
numerical approximations. However, \cite{CW,HMY,WMS} only dealt with
the case that the coefficient of the diffusion part is at most
linear growth, that is, there exists $K>0$ such that
\begin{equation}\label{linear}
|g(x)|\le K|x|, \forall x\in \mathbb{R}^d.
\end{equation}

This condition excludes the case when the coefficient $g$ is
super-linearly growing (that is, $g(x)=C|x|^\gamma,\ \gamma>1)$. In
Mao and Szpruch \cite{MS}, authors examined the globally almost sure
asymptotic stability of the $\theta$-EM scheme (\ref{SEM1}), they
presented a rather weak sufficient condition to ensure that the
$\theta$-EM solution is almost surely stable when
$\frac{1}{2}<\theta\le 1$, but they didn't give the convergence rate
of the solution to zero explicitly. In \cite{ZW}, the authors
studied the mean square exponential stability of $\theta$-EM scheme
systematically, they proved that if $0\le\theta<\frac{1}{2},$ the
$\theta$-EM scheme preserves mean square exponential stability under
the linear growth condition for both the drift term and the
diffusion term; if $\frac{1}{2}<\theta\le 1,$ the $\theta$-EM
preserves mean square exponential stability without the linear
growth condition for the drift term (the linear growth condition for
the diffusion term is still necessary), exponential stability for
the case $\theta=\frac{1}{2}$ is not studied there.

However, to the best of our knowledge, there are few results devoted
to the exponential stability of the numerical solutions when the
coefficient of the diffusion term does not satisfy the linear growth
condition, which is one of the main motivations of this work.

Recently, in \cite{LFM}, Liu et al examined the polynomial stability
of numerical solutions of SDEs (\ref{sde}). They considered the
polynomial stability of both the classical and backward
Euler-Maruyama approximation. The condition on diffusion coefficient
$g$ is bounded with respect to variable $x$. This condition excludes
the case that $g$ is unbounded with respect to variable $x$. It
immediately raises the question of whether we can relax this
condition. This is the other main motivation of this work.

To study the polynomial stability of equation (\ref{SEM}), we consider the
following condition:
\begin{equation}\label{c1}
2\langle x,f(x,t)\rangle+|g(x,t)|^2\le C(1+t)^{-K_1}- K_1(1+t)^{-1}|x|^2,
\forall t\ge0, x\in\mathbb{R}^d,
\end{equation}
where $K_1,  C$ are positive constants, and $K_1>1,$ $\langle
\cdot,\cdot\rangle$ stands for the inner product in $\mathbb{R}^d$
and  $|\cdot|$ denotes the both the Euclidean vector norm and the
Hilbert-Schmidt matrix norm.

To study the exponential stability of equation (\ref{SEM}),  we need stronger
condition on the coefficients,
\begin{equation}\label{c2}
2\langle x,f(x,t)\rangle+|g(x,t)|^2\le -C|x|^2, \forall  x\in\mathbb{R}^d,
\end{equation}
where $C>0$ is a constant.

Define an operator $L$ by
$$\aligned LV(x,t):&=\frac{\partial}{\partial t}V(x,t)+\sum_{i=1}^df^i(x,t)
\frac{\partial}{\partial x_i}V(x,t)\\&
\quad+\frac{1}{2}\sum_{i,j=1}^d\sum_{k=1}^m g^{ik}(x,t)g^{jk}(x,t)\frac{\partial^2}
{\partial x_i\partial x_j}V(x,t),\endaligned$$
where $V(x,t): \mathbb{R}^d\times\mathbb{R}^+\rightarrow\mathbb{R}^+$ has continuous
second-order partial derivatives in $x$ and first-order partial derivatives in $t.$

It is clear that under condition (\ref{local}) and (\ref{c1}) (or
(\ref{c2})), there exists a unique global solution of equation
(\ref{sde}). By taking $V(x,t)=(1+t)^m|x|^2,$ or $V(x,t)=|x|^2,$
respectively, it is easy to see that under condition (\ref{c1}) the
true solution $X_t(x_0)$ of equation (\ref{sde}) is mean square
polynomially stable (see Liu and Chen \cite{LC} Theorem 1.1) or mean
square exponentially stable under condition (\ref{c2}) (the proof is
the same as Higham et al, see \cite{HMY} Appendix A). So a natural
question raises: Whether $\theta$-EM method can reproduce the
polynomial and exponential stability of the solution of (\ref{sde}).

If $\frac{1}{2} <\theta\le 1$, we will study the polynomial
stability and exponential stability of $\theta$-EM scheme
(\ref{SEM}) under conditions (\ref{c1}) and (\ref{c2}) respectively.
For the exponential stability, we first investigate the mean square
exponential stability, then we derive the almost sure exponential
stability by Borel-Cantelli lemma.

If $0\le\theta\le\frac{1}{2},$ besides condition (\ref{c1})
(respectively, (\ref{c2})), linear growth condition for the drift
term is also needed to ensure the corresponding stability, that is,
there exists $K>0$ such that
\begin{equation}\label{growth}
|f(x,t)|\le K(1+t)^{-\frac{1}{2}}|x|
\end{equation}
for polynomial stability case and
\begin{equation}\label{growth1}
|f(x,t)|\le K|x|, \forall x\in \mathbb{R}^d
\end{equation}
for exponential stability case. Notice that condition (\ref{growth})
is strictly weaker than condition (2.4) in \cite{LFM}.

The main feature of this paper is that we consider conditions in
which both diffusion and drift coefficients are involved, which give
weaker sufficient conditions than known ones, while in most of the
preceding studies, such conditions have been provided as separate
ones for diffusion coefficients and drift coefficients.

The rest of the paper is organized as follows. In Section 2, we give
some lemmas which will be used in the following sections to prove
the stability results. In Section 3 we study the polynomial
stability of the $\theta$-EM scheme. Our method hinges on various
properties of the gamma function and the ratios of gamma functions.
We show that when $\frac{1}{2}<\theta\le 1$, the polynomial
stability of the $\theta$-EM scheme holds under condition (\ref{c1})
plus one sided Lipschitz condition on $f$; when
$0\le\theta\le\frac{1}{2},$ the linear growth condition for the
drift term $f$ is also needed. In Section 4, we investigate the
exponential stability of the $\theta$-EM scheme for all
$0\le\theta\le 1$. Finally, we give in Section 5 some non stability
results and counter examples to support our conclusions.

\section{Preliminary}

To ensure that  the semi-implicit $\theta$-EM scheme is well
defined, we need the first two lemmas.The first lemma gives the
existence and uniqueness of the solution of the equation $F(x)=b.$
We can prove the existence and uniqueness of the solution of the
$\theta$-EM scheme based on this lemma.

\begin{Lemma}\label{l1}
Let $F$ be the vector field on $\mathbb{R}^d$ and consider the equation
\begin{equation}\label{e1}
F(x)=b
\end{equation}
for a given $b\in\mathbb{R}^d$. If $F$ is monotone, that is,
$$\langle x-y, F(x)-F(y)\rangle>0$$
for all $x,y\in\mathbb{R}^d,x\neq y$,  and $F$ is continuous and coercive, that is,
$$\lim_{|x|\rightarrow\infty}\frac{\langle x, F(x)\rangle}{|x|}=\infty,$$
then for every $b\in\mathbb{R}^d,$ equation (\ref{e1}) has a unique solution $x\in\mathbb{R}^d$.
\end{Lemma}

This lemma follows directly from Theorem 26.A in \cite{Zeidler}.

Consider the following one sided Lipschitz condition on $f$: There exists $L>0$ such that

\begin{equation}\label{c3}
\langle x-y, f(x,t)-f(y,t)\rangle\le L|x-y|^2.
\end{equation}

\begin{Lemma}\label{l2}
Define \begin{equation}\label{dingyi}F(x,t):=x-\theta\Delta t
f(x,t), \forall t>0, x\in\mathbb{R}^d.\end{equation} Assume
conditions (\ref{c1}) and (\ref{c3}) and $\Delta t$ is small enough
such that $\Delta t<\frac{1}{\theta L}$. Then for any $t>0$ and
$b\in\mathbb{R}^d,$ there is a unique solution of equation
$F(x,t)=b.$
\end{Lemma}

By this Lemma, we know that the $\theta$-EM scheme is well defined
under conditions (\ref{c1}) and (\ref{c3}) for $\Delta t$  small
enough.

The proof of Lemma \ref{l2} is the same as that of Lemma 3.4 in
\cite{LFM} and Lemma 3.3 in \cite{MS2}, just notice that condition
(\ref{c3}) implies $\langle x-y,F(x,t)-F(y,t)\rangle>0$, and
(\ref{c1}) (or (\ref{c2})) implies $\langle
x,F(x,t)\rangle\rightarrow\infty$ as $x\rightarrow\infty$. Notice
also that our condition (\ref{c3}) is weaker than (2.3) in
\cite{LFM}.

We also need the following two lemmas to study the polynomial stability
of the $\theta$-EM scheme.

\begin{Lemma}\label{l3}
Given $\alpha>0$ and $\beta\ge 0,$ if there exists a $\delta$ such that
$0<\delta<\alpha^{-1}$, then
$$\prod_{i=a}^b\Big(1-\frac{\alpha\delta}{1+(i+\beta)\delta}\Big)=
\frac{\Gamma(b+1+\delta^{-1}+\beta-\alpha)}{\Gamma(b+1+\delta^{-1}+\beta)}
\times\frac{\Gamma(a+\delta^{-1}+\beta)}{\Gamma(a+\delta^{-1}+\beta-\alpha)},$$
where $0\le a\le b,$ $\Gamma(x):=\int_0^\infty y^{x-1}e^{-y}dy.$
\end{Lemma}

\begin{Lemma}\label{l4}
For any $x>0,$ if $0<\eta<1,$ then
$$\frac{\Gamma(x+\eta)}{\Gamma(x)}<x^\eta,$$
and if $\eta>1,$ then
$$\frac{\Gamma(x+\eta)}{\Gamma(x)}>x^\eta.$$

\end{Lemma}

The proof of Lemmas \ref{l3} and \ref{l4} could be found in \cite{LFM}.

\section{Polynomial stability of $\theta$-EM solution (\ref{SEM})}

We are now in the position to give the polynomial stability of $\theta$-EM solution
(\ref{SEM}). First, we consider the case $\frac{1}{2}<\theta\le 1.$ We have the following

\begin{Theorem}\label{polynomial}
Assume that conditions (\ref{c1}) and (\ref{c3}) hold. If
$\frac{1}{2}<\theta\le 1,$ then for any $0<\varepsilon<K_1-1$, we
can choose $\Delta t$ small enough such that the $\theta$-EM
solution satisfies
\begin{equation}\label{bu}
\limsup_{k\rightarrow \infty}\frac{\log \mathbb{E}|X_k|^2}{\log
k\Delta t}\le -(K_1-1-\varepsilon)
\end{equation}
for any initial value $X_0=x_0\in \mathbb{R}^d.$
\end{Theorem}

\textbf{Proof} We first prove that condition (\ref{c1}) implies that
for $\Delta t$ small enough,
\begin{equation}\label{ineq}
2\langle x,f(x,t)\rangle+|g(x,t)|^2+(1-2\theta)\Delta t|f(x,t)|^2\le
C(1+t)^{-K_1}-(K_1-\varepsilon)(1+t)^{-1}|F(x,t)|^2 \end{equation}
holds for $\forall t\ge0, x\in\mathbb{R}^d.$ Here and in the
following, $F$ is defined by (\ref{dingyi}).

In fact, we only need to show that
$$(2\theta-1)\Delta t|f(x,t)|^2- (K_1-\varepsilon)(1+t)^{-1}|F(x,t)|^2\ge - K_1(1+t)^{-1}|x|^2.$$

On the other hand, by the definition of $F(x,t),$ we have
$$\aligned &\quad\ (2\theta-1)\Delta t|f(x,t)|^2- (K_1-\varepsilon)(1+t)^{-1}|F(x,t)|^2\\&
=(2\theta-1)\Delta t|f(x,t)|^2-
(K_1-\varepsilon)(1+t)^{-1}[|x|^2-2\theta\Delta t\langle x,
f(x,t)\rangle+\theta^2\Delta t^2|f(x,t)|^2]\\& =[(2\theta-1)\Delta
t-(K_1-\varepsilon)(1+t)^{-1}\theta^2\Delta
t^2]|f(x,t)|^2\\&\quad+2(K_1-\varepsilon)(1+t)^{-1}\theta\Delta
t\langle x,f(x,t)\rangle- (K_1-\varepsilon)(1+t)^{-1}|x|^2\\&
=a|f(x,t)+bx|^2-(ab^2+(K_1-\varepsilon)(1+t)^{-1})|x|^2,\endaligned
$$
where
$$a:=(2\theta-1)\Delta t-(K_1-\varepsilon)(1+t)^{-1}\theta^2\Delta t^2,\quad
b:=\frac{(K_1-\varepsilon)(1+t)^{-1}\theta\Delta t}{a}.$$

Since
$$a\ge (2\theta-1)\Delta t-(K_1-\varepsilon)\theta^2\Delta t^2
=\Delta t(2\theta-1-(K_1-\varepsilon)\theta^2\Delta t),$$ we can
choose $\Delta t$ small enough (for example $\Delta t\le
\min\{\frac{1}{\theta L},\frac{(2\theta-1)(\varepsilon\wedge
1)}{K_1(K_1-\varepsilon)\theta^2}\}$) such that $a\ge 0$ and
$ab^2\le \frac{\varepsilon}{1+t}.$ Then we have

$$\aligned \ (2\theta-1)\Delta t|f(x,t)|^2- (K_1-\varepsilon)(1+t)^{-1}|F(x,t)|^2&
\ge -(\frac{b^2}{a}+(K_1-\varepsilon)(1+t)^{-1})|x|^2\\&
\ge -K_1(1+t)^{-1}|x|^2.\endaligned
$$

So we complete the proof of inequality (\ref{ineq}).

Now by the definition of $F(x,t)$, it follows that
$$F(X_{k+1},(k+1)\Delta t)=F(X_{k},k\Delta t)+f(X_k,k\Delta t)\Delta t
+g(X_k,k\Delta t)\Delta B_k.$$

So we have
\begin{equation}\label{F}\aligned |F(X_{k+1},(k+1)\Delta t)|^2&=
[2\langle X_k,f(X_k,k\Delta t)\rangle+|g(X_{k},k\Delta t)|^2+(1-2\theta)
|f(X_{k},k\Delta t)|^2\Delta t]\Delta t\\&\quad+|F(X_{k},k\Delta t)|^2+M_k,\endaligned\end{equation}
where
\begin{equation}\label{M}\aligned M_k&:=|g(X_{k},k\Delta t)\Delta B_k|^2
-|g(X_{k},k\Delta t)|^2\Delta t+2\langle F(X_{k},k\Delta t),g(X_{k},k\Delta t)\Delta B_k\rangle\\&
\quad+2\langle f(X_{k},k\Delta t)\Delta t,g(X_{k},k\Delta t)\Delta B_k\rangle.\endaligned\end{equation}

Notice that
$$\mathbb{E}(M_k|\mathscr{F}_{k\Delta t})=0.$$

Then by condition (\ref{c1}) and inequality (\ref{ineq}), we have
$$ \mathbb{E}(|F(X_{k+1},(k+1)\Delta t)|^2|\mathscr{F}_{k\Delta t})\le
(1-\frac{(K_1-\varepsilon)\Delta t}{1+k\Delta t})|F(X_{k},k\Delta t)|^2
+\frac{C\Delta t}{(1+k\Delta t)^{K_1-\varepsilon}}.$$

We can get by iteration that
$$\aligned\mathbb{E}(|F(X_{k},k\Delta t)|^2)&\le \Big(\prod_{i=0}^{k-1}
(1-\frac{(K_1-\varepsilon)\Delta t}{1+k\Delta t})\Big)|F(x_0,0)|^2\\&\quad
+\sum_{r=0}^{k-1}\Big(\prod_{i=r+1}^{k-1}(1-\frac{(K_1-\varepsilon)\Delta t}{1+i\Delta t})\Big)
\frac{C\Delta t}{(1+r\Delta t)^{K_1-\varepsilon}}.\endaligned$$

Then by Lemma \ref{l3},
\begin{equation}\label{budeng}\aligned\mathbb{E}(|F(X_{k},k\Delta t)|^2)&\le
\frac{\Gamma(k+\frac{1}{\Delta t}-(K_1-\varepsilon))\Gamma(\frac{1}{\Delta t})}{\Gamma(k+\frac{1}{\Delta t})
\Gamma(\frac{1}{\Delta t}-(K_1-\varepsilon))}|F(x_0,0)|^2\\&\quad
+C\Delta t\sum_{r=0}^{k-1}\frac{\Gamma(k+\frac{1}{\Delta t}-(K_1-\varepsilon))
\Gamma(r+1+\frac{1}{\Delta t})}{\Gamma(k+\frac{1}{\Delta t})\Gamma(r+1+\frac{1}{\Delta t}
-(K_1-\varepsilon))}(1+r\Delta t)^{-(K_1-\varepsilon)}.\endaligned\end{equation}

On the other hand, since $K_1-\varepsilon>1,$ by Lemma \ref{l4} or
\cite{LFM} one can see that
\begin{equation}\label{kongzhi1}\frac{\Gamma(k+\frac{1}{\Delta t}-(K_1-\varepsilon))
\Gamma(\frac{1}{\Delta t})}{\Gamma(k+\frac{1}{\Delta t})\Gamma(\frac{1}{\Delta t}-
(K_1-\varepsilon))}\le ((k-(K_1-\varepsilon))\Delta t+1)^{-(K_1-\varepsilon)}\end{equation}
and that
\begin{equation}\label{kongzhi2}\frac{\Gamma(k+\frac{1}{\Delta t}-(K_1-\varepsilon))
\Gamma(r+1+\frac{1}{\Delta t})}{\Gamma(k+\frac{1}{\Delta t})\Gamma(r+1+\frac{1}{\Delta t}-
(K_1-\varepsilon))}\le ((k-(K_1-\varepsilon))\Delta t+1)^{-(K_1-\varepsilon)}((r+1)\Delta t+1)^{K_1-\varepsilon}.\end{equation}

Substituting (\ref{kongzhi1}) and (\ref{kongzhi2}) into inequality (\ref{budeng}) yields
\begin{equation}\aligned\mathbb{E}(|F(X_{k},k\Delta t)|^2)&\le((k-(K_1-\varepsilon))\Delta t+1)^{-(K_1-\varepsilon)}
|F(x_0,0)|^2\\&\quad+C\Delta
t\sum_{r=0}^{k-1}((k-(K_1-\varepsilon))\Delta
t+1)^{-(K_1-\varepsilon)} \frac{((r+1)\Delta
t+1)^{K_1-\varepsilon}}{(1+r\Delta t)^{K_1-\varepsilon}}\\&
\le2^{K_1-\varepsilon}(k\Delta
t+1)^{-(K_1-\varepsilon)}\Big[|F(x_0,0)|^2+C\Delta t\sum_{r=0}^{k-1}
\Big(\frac{(r+1)\Delta t+1}{1+r\Delta
t}\Big)^{K_1-\varepsilon}\Big]\\&\le2^{K_1-\varepsilon}(k\Delta
t+1)^{-(K_1-\varepsilon)} [|F(x_0,0)|^2+C\cdot
2^{K_1-\varepsilon}k\Delta t]\\&\le
2^{K_1-\varepsilon}(|F(x_0,0)|^2+C\cdot 2^{K_1-\varepsilon})
(k\Delta t+1)^{-(K_1-\varepsilon)+1}.\endaligned\end{equation} We
have used the fact that $(k-(K_1-\varepsilon))\Delta t+1\ge
\frac{1}{2}(k\Delta t+1) $ for small $\Delta t$ in second inequality
and that $((r+1)\Delta t+1)/(1+r\Delta t)\le 2$ in the third
inequality.

Now by condition (\ref{c1}),
$$\aligned |F(x,t)|^2&=|x|^2-2\theta\Delta t\langle x,f(x,t)\rangle+\theta^2\Delta t^2|f(x,t)|^2\\&
\ge |x|^2- C(1+t)^{-K_1}\theta\Delta t+ K_1(1+t)^{-1}|x|^2\theta\Delta t+\theta^2\Delta t^2|f(x,t)|^2\\&
\ge |x|^2- C(1+t)^{-K_1}\theta\Delta t\ge |x|^2- C(1+t)^{-(K_1-\varepsilon)}\theta\Delta t.\endaligned$$

Therefore, for small enough $\Delta t,$
$$\aligned\mathbb{E}(|X_{k}|^2)&\le\mathbb{E}(|F(X_{k},k\Delta t)|^2)+C(1+k\Delta t)^{-(K_1-\varepsilon)}\theta\Delta t\\&
\le 2^{K_1-\varepsilon}(|F(x_0,0)|^2+C\cdot
2^{K_1-\varepsilon})(k\Delta t+1)^{-(K_1-\varepsilon)+1}+C(1+k\Delta
t)^{-(K_1-\varepsilon)}\theta\Delta t\\& \le
2^{K_1-\varepsilon}(|F(x_0,0)|^2+C\cdot2^{K_1-\varepsilon}+C\theta)(k\Delta
t+1)^{-(K_1-\varepsilon)+1}.\endaligned$$ Namely, the $\theta$-EM
solution of (\ref{sde}) is mean square polynomial stable with rate
no greater than $-(K_1-1-\varepsilon)$ when $\frac{1}{2}<\theta\le
1$ and $\Delta t$ is small enough.

We complete the proof. $\square$

\begin{Remark}
Notice that we can not let $\varepsilon\rightarrow0$ in (\ref{bu})
since $\Delta t$ depends on $\varepsilon.$ Moreover, our condition
(\ref{c1}) could cover conditions (2.5) and (2.6) (even though not
entirely. They need $K_1>0.5$, but our $K_1>1$) for the polynomial
stability of backward EM approximation of SDE (\ref{sde}).
\end{Remark}

Now let us consider the case $0\le\theta\le\frac{1}{2}.$ We have

\begin{Theorem}\label{polynomial1}
Assume that conditions (\ref{c1}), (\ref{growth}) and (\ref{c3})
hold. If $0\le\theta\le\frac{1}{2},$ then for any
$0<\varepsilon<K_1-1$, we can choose $\Delta t$ small enough such
that the $\theta$-EM solution satisfies
\begin{equation}
\limsup_{k\rightarrow \infty}\frac{\log \mathbb{E}|X_k|^2}{\log
k\Delta t}\le -(K_1-1-\varepsilon)
\end{equation}
for any initial value $X_0=x_0\in \mathbb{R}^d.$
\end{Theorem}

\textbf{Proof}
Notice that in this case
$$\aligned &\quad\ (2\theta-1)\Delta t|f(x,t)|^2- (K_1-\varepsilon)(1+t)^{-1}|F(x,t)|^2\\&
=(2\theta-1)\Delta t|f(x,t)|^2-
(K_1-\varepsilon)(1+t)^{-1}[|x|^2-2\theta\Delta t\langle
x,f(x,t)\rangle +\theta^2\Delta t^2|f(x,t)|^2]\\&
=[(2\theta-1)\Delta t-(K_1-\varepsilon)(1+t)^{-1}\theta^2\Delta
t^2]|f(x,t)|^2\\&\quad+2 (K_1-\varepsilon) (1+t)^{-1}\theta\Delta
t\langle x,f(x,t)\rangle- (K_1-\varepsilon)(1+t)^{-1}|x|^2\\& \ge
aK^2(1+t)^{-1}|x|^2-2
K(K_1-\varepsilon)(1+t)^{-\frac{3}{2}}\theta\Delta
t|x|^2-(K_1-\varepsilon)(1+t)^{-1}|x|^2,\endaligned
$$
where
$$a:=(2\theta-1)\Delta t-(K_1-\varepsilon)(1+t)^{-1}\theta^2\Delta t^2\le 0.$$

Thus, we can choose $\Delta t$ small enough such that
$$aK^2(1+t)^{-1}-2 KK_1(1+t)^{-\frac{3}{2}}\theta\Delta t-(K_1-\varepsilon)(1+t)^{-1}\ge -K_1(1+t)^{-1}.$$

Therefore, by condition (\ref{c1}), we have
$$ \mathbb{E}(|F(X_{k+1},(k+1)\Delta t)|^2|\mathscr{F}_{k\Delta t})\le (1-\frac{(K_1-\varepsilon)
\Delta t}{1+k\Delta t})|F(X_{k},k\Delta t)|^2+\frac{C\Delta t}{(1+k\Delta t)^{K_1-\varepsilon}}.$$

Then by the same argumentation as Theorem \ref{polynomial}, we have
$$\aligned |F(x,t)|^2&=|x|^2-2\theta\Delta t\langle x,f(x,t)\rangle+\theta^2\Delta t^2|f(x,t)|^2\\&
\ge |x|^2- C(1+t)^{-K_1}\theta\Delta t+ K_1(1+t)^{-1}|x|^2\theta\Delta t\\&
\ge |x|^2- C(1+t)^{-K_1}\theta\Delta t\ge |x|^2- C(1+t)^{-(K_1-\varepsilon)}\theta\Delta t.\endaligned$$

Therefore, for small enough $\Delta t,$ we can derive in the same way as in proof of Theorem \ref{polynomial} that
$$\aligned\mathbb{E}(|X_{k}|^2)&\le\mathbb{E}(|F(X_{k},k\Delta t)|^2)+C(1+k\Delta t)^{-(K_1-\varepsilon)}\theta\Delta t\\&
\le 2^{K_1-\varepsilon}(|F(x_0,0)|^2+C\cdot
2^{K_1-\varepsilon})(k\Delta t+1)^{-(K_1-\varepsilon)+1}+C(1+k\Delta
t)^{-(K_1-\varepsilon)}\theta\Delta t\\& \le
2^{K_1-\varepsilon}(|F(x_0,0)|^2+C\cdot2^{K_1-\varepsilon}+C\theta)(k\Delta
t+1)^{-(K_1-\varepsilon)+1}.\endaligned$$

Namely, the $\theta$-EM solution of (\ref{sde}) is mean square
polynomial stable with rate no greater than $-(K_1-1-\varepsilon)$
when $0\le\theta\le\frac{1}{2}$ and $\Delta t$ is small enough.

We complete the proof. $\square$

\begin{Remark}
In \cite{LFM} Condition 2.3, authors gave the sufficient conditions
on coefficients $f$ and $g$ separately for the polynomial stability
of the classical EM scheme, their conditions (2.5) and (2.6) hold
for $K_1>1$ and $C>0,$ then it is easy to see that our condition
(\ref{c1}) holds automatically for the same $K_1$ and $C,$ and our
condition (\ref{growth}) is strictly weaker than (2.4). Therefore,
we have improved Liu et al and generalized it to
$0\le\theta\le\frac{1}{2}.$
\end{Remark}

\section{Exponential stability of $\theta$-EM solution (\ref{SEM})}

Now let us consider the exponential stability of $\theta$-EM solution of  (\ref{sde}).
When SDE (\ref{sde}) goes back to time homogeneous case, that is,
\begin{equation}\label{sde1}dX_t=f(X_t)dt+g(X_t)dB_t,\
X_0=x_0\in \mathbb{R}^d,a.s. \end{equation}

The corresponding $\theta$-EM approximation becomes to
\begin{equation}\label{SEM1} X_{k+1}=X_k+[(1-\theta)f(X_k)+\theta f(X_{k+1})]\Delta t
+g(X_k)\Delta B_k.\end{equation}

In \cite{MS}, Mao and Szpruch gave a sufficient condition ensuring
that the almost sure stability of $\theta$-EM solution of
(\ref{sde1}) holds in the case that $\frac{1}{2}<\theta\le 1$.
However they didn't reveal the rate of convergence. Their method of
the proof is mainly based on the discrete semi martingale
convergence theorem. We will study the exponential stability
systematically for $0\le\theta\le 1$ for the time inhomogeneous
case. We first prove the mean square exponential stability, then we
prove the almost sure stability by Borel-Cantelli lemma.

\begin{Theorem}\label{exponential}
Assume that conditions (\ref{c2}) and (\ref{c3}) hold. Then for any
$\frac{1}{2}<\theta\le 1$ and $0<\varepsilon<1$, we can choose
$\Delta t$ small enough such that the $\theta$-EM solution satisfies
\begin{equation}
\limsup_{k\rightarrow \infty}\frac{\log \mathbb{E}|X_k|^2}{k\Delta t}\le -C(1-\varepsilon)
\end{equation}
for any initial value $X_0=x_0\in \mathbb{R}^d$ and
\begin{equation}\label{ex}
\limsup_{k\rightarrow \infty}\frac{\log |X_k|}{k\Delta t}\le -\frac{C(1-\varepsilon)}{2}\quad a.s.
\end{equation}

\end{Theorem}

\textbf{Proof} Define $F(x,t)$ as in Lemma \ref{l2}. We have
$$\aligned &\quad\ (2\theta-1)\Delta t|f(x,t)|^2-C(1-\varepsilon)|F(x,t)|^2\\&
=(2\theta-1)\Delta t|f(x,t)|^2- C(1-\varepsilon)[|x|^2-2\theta\Delta t\langle x,f(x,t)\rangle
+\theta^2\Delta t^2|f(x,t)|^2]\\&
=[(2\theta-1)\Delta t-C\theta^2\Delta t^2(1-\varepsilon)]|f(x,t)|^2+2C\theta\Delta t
(1-\varepsilon)\langle x,f(x,t)\rangle-C(1-\varepsilon)|x|^2\\&
=a|f(x,t)+bx|^2-(ab^2+C(1-\varepsilon))|x|^2,\endaligned
$$
where
$$a:=(2\theta-1)\Delta t-C\theta^2\Delta t^2(1-\varepsilon), \quad
b:=\frac{C\theta\Delta t(1-\varepsilon)}{a}.$$

We can choose $\Delta t$ small enough (for example $\Delta t\le \min
\{\frac{1}{\theta L},\frac{\varepsilon(2\theta-1)}{C(1-\varepsilon)\theta^2}\}$)
such that $a\ge 0$ and $ab^2\le C\varepsilon$, and therefore
$$(2\theta-1)\Delta t|f(x,t)|^2-C(1-\varepsilon)|F(x,t)|^2\ge -C|x|^2.$$

Then by condition (\ref{c2}), we can prove that
\begin{equation}\label{ineq1}
2\langle x,f(x,t)\rangle+|g(x,t)|^2+(1-2\theta)\Delta t|f(x,t)|^2\le
-C(1-\varepsilon)|F(x,t)|^2 \end{equation}
holds for $\forall x\in\mathbb{R}^d.$

Therefore, by  (\ref{F}), for small enough $\Delta t$ ($\Delta t\le
\frac{1}{\theta L}\wedge\frac{\varepsilon(2\theta-1)}{C\theta^2(1-\varepsilon)}
\wedge\frac{1}{C(1-\varepsilon)}$),
$$\aligned \mathbb{E}(|F(X_{k+1},(k+1)\Delta t)|^2)
\le \mathbb{E}(|F(X_{k},k\Delta t)|^2)(1-C(1-\varepsilon)\Delta t).\endaligned$$

So we have
\begin{equation}\label{mean1}\aligned \mathbb{E}(|X_{k}|^2)\le
\mathbb{E}(|F(X_{k},k\Delta t)|^2) \le
|F(x_0,0)|^2(1-C(1-\varepsilon)\Delta t)^k.\endaligned\end{equation}
or
\begin{equation}\label{mean}\aligned
\mathbb{E}(|X_{k}|^2)\le |F(x_0,0)|^2e^{-C(1-\varepsilon)k\Delta t},
\forall k\ge 1.\endaligned\end{equation}

The first inequality of (\ref{mean1}) holds because of condition (\ref{c2}).
Thus, the $\theta$-EM solution of (\ref{sde1}) is mean square exponential
stable when $\frac{1}{2}<\theta\le 1$ and $\Delta t$ is small enough.

On the other hand, by Chebyshev inequality, inequality (\ref{mean}) implies that
$$P(|X_k|^2>k^2e^{-kC(1-\varepsilon)\Delta t})\le\frac{|F(x_0,0)|^2}{k^2}, \forall k\ge 1.$$

Then by Borel-Cantelli lemma, we see that for almost all $\omega\in\Omega$
\begin{equation}\label{bound}|X_k|^2\le k^2e^{-kC(1-\varepsilon)\Delta t}\end{equation}
holds for all but finitely many $k$. Thus, there exists a $k_0(\omega),$ for
all $\omega\in\Omega$ excluding a $P$-null set, for which (\ref{bound}) holds whenever $k\ge k_0$.

Therefore, for almost all $\omega\in\Omega$,
\begin{equation}\label{bound1}\frac{1}{k\Delta t}\log|X_k|\le-\frac{C(1-\varepsilon)}{2}
+\frac{\log k}{k\Delta t}\end{equation}
whenever $k\ge k_0$. Letting $k\rightarrow\infty$ we obtain (\ref{ex}).

The proof is then complete. $\square$

If $0\le\theta\le\frac{1}{2},$ then we have the following

\begin{Theorem}\label{exponential1}
Assume that conditions (\ref{c2}), (\ref{growth1}) and (\ref{c3}) hold. Then for any
$0<\varepsilon<1$, we can choose $\Delta t$ small enough
such that the $\theta$-EM solution satisfies
\begin{equation}
\limsup_{k\rightarrow \infty}\frac{\log \mathbb{E}|X_k|^2}{k\Delta t}\le -C(1-\varepsilon)
\end{equation}
for any initial value $X_0=x_0\in \mathbb{R}^d$ and
\begin{equation}\label{ex1}
\limsup_{k\rightarrow \infty}\frac{\log |X_k|}{k\Delta t}\le -\frac{C(1-\varepsilon)}{2}\quad a.s.
\end{equation}

\end{Theorem}

\textbf{Proof} By  the same argument as Theorem \ref{polynomial1}, we have
$$\aligned &\quad\ (2\theta-1)\Delta t|f(x,t)|^2-C(1-\varepsilon)|F(x,t)|^2\\&
=(2\theta-1)\Delta t|f(x,t)|^2- C(1-\varepsilon)[|x|^2-2\theta\Delta t\langle x,f(x,t)\rangle
+\theta^2\Delta t^2|f(x,t)|^2]\\&
=[(2\theta-1)\Delta t-C\theta^2\Delta t^2(1-\varepsilon)]|f(x,t)|^2+2C\theta\Delta t
(1-\varepsilon)\langle x,f(x,t)\rangle-C(1-\varepsilon)|x|^2\\&
\ge aK^2|x|^2-2KC\theta\Delta t(1-\varepsilon)|x|^2-C(1-\varepsilon)|x|^2\endaligned
$$
since
$$a:=(2\theta-1)\Delta t-C\theta^2\Delta t^2(1-\varepsilon)\le 0.$$
We have used condition (\ref{growth}) in the last inequality.

We can choose $\Delta t$ small enough such that
$$\Delta t\le \frac{1}{\theta L}\wedge\frac{K(1-2\theta)+2C\theta(1-\varepsilon)}{KC(1-\varepsilon)\theta^2}
\wedge\frac{C\varepsilon}{2(K^2(1-2\theta)+2KC\theta(1-\varepsilon))},$$ and thus
$$aK^2-2KC\theta\Delta t(1-\varepsilon)\ge -C\varepsilon.$$

Then we have
$$(2\theta-1)\Delta t|f(x,t)|^2-C(1-\varepsilon)|F(x,t)|^2\ge -C|x|^2.$$

Then by condition (\ref{c2}), we can prove that
\begin{equation}\label{ineq1}
2\langle x,f(x,t)\rangle+|g(x,t)|^2+(1-2\theta)\Delta t|f(x,t)|^2\le -C(1-\varepsilon)|F(x,t)|^2 \end{equation}
holds for $\forall x\in\mathbb{R}^d.$

Therefore, for small enough $\Delta t$,
$$\aligned \mathbb{E}(|F(X_{k+1},(k+1)\Delta t)|^2)
\le \mathbb{E}(|F(X_{k},k\Delta t)|^2)(1-C(1-\varepsilon)\Delta t).\endaligned$$

So we have
\begin{equation}\label{mean2}\aligned \mathbb{E}(|X_{k}|^2)\le \mathbb{E}(|F(X_{k},k\Delta t)|^2)
\le |F(x_0,0)|^2(1-C(1-\varepsilon)\Delta t)^k.\endaligned\end{equation}
or
\begin{equation}\label{mean3}\aligned
\mathbb{E}(|X_{k}|^2)\le |F(x_0,0)|^2e^{-C(1-\varepsilon)k\Delta t}, \forall k\ge 1.\endaligned\end{equation}

The first inequality of (\ref{mean1}) holds because of condition (\ref{c2}). Thus, the $\theta$-EM
solution of (\ref{sde1}) is mean square exponential stable when $\frac{1}{2}<\theta\le 1$ and $\Delta t$
is small enough.

From (\ref{mean}) we can show the almost sure stability assertion (\ref{ex1}) in the same
way as in the proof of Theorem \ref{exponential}.

The proof is complete. $\square$

\section{Non stability results and counter examples}

In this section we will give some non stability results for the classical EM scheme
and counter examples to support our conclusions. We show that there are cases that
our assertion works while the assertions in the literature do not work.

Let us consider the following 1-dimensional stochastic differential equations:
\begin{equation}\label{SDE} dX_t=(aX_t+b|X_t|^{q-1}X_t)dt+c|X_t|^\gamma dB_t.\ X_0=x_0(\neq0)\in \mathbb{R}.\end{equation}

When $b\le0, q> 0$ and $\gamma\ge\frac{1}{2},$ by Gy\"ongy and
Krylov \cite{GK} Corollary 2.7 (see also \cite{IW,Lan,Lan1,RY}),
there is a unique global solution of equation (\ref{SDE}). Here
$|x|^{q-1}x:=0$ if $x=0$. For this equation, if $q=2\gamma-1, a<0$
and $2b+c^2\le 0,$ then condition (\ref{c2}) is automatically
satisfied. Therefore, the true solution of SDE (\ref{SDE}) is mean
square exponentially stable.

Now let us consider the corresponding Euler-Maruyama approximation:
\begin{equation}\label{EM} X_{k+1}=X_k+(aX_k+b|X_k|^{q-1}X_k)\Delta t+c|X_k|^\gamma \Delta B_k.\end{equation}

For the classical EM approximation $X_k$, we have the following

\begin{Lemma}\label{divergence1}
Suppose $q>1, q>\gamma.$ If $\Delta t>0$ is small enough, and
$$|X_1|\ge\frac{2^\frac{q+2}{q-1}}{(|b|\Delta t)^\frac{1}{q-1}},$$
then for any $K\ge 1,$ there exists a positive number $\alpha$ such that
\[P(|X_k|\ge\frac{2^{k+\frac{3}{q-1}}}{(|b|\Delta t)^\frac{1}{q-1}},\forall
1\le k\le K)\ge\exp(-4e^{-\alpha/\sqrt{\Delta t}})>0,\]
where $\alpha:=\frac{2^{\frac{(q-\gamma)(q+2)}{q-1}}}{2|c|}(1\wedge((q-\gamma)\log 2)).$
\end{Lemma}
That is, no matter what values $a,b,c$ take, by taking the initial
value and the step size suitably, the numerical approximation
solution of SDE (\ref{SDE}) is divergent with a positive probability
when $q>1, q>\gamma.$

\textbf{Proof of Lemma \ref{divergence1}}: According to (\ref{EM}),
\[\aligned|X_{k+1}|&=|X_k|\Big| b|X_k|^{q-1}\Delta t+c\cdot \textrm{sgn}
(X_k)|X_k|^{\gamma-1}\Delta B_k+1+a\Delta t\Big|\\&
\ge |X_k|\Big||b||X_k|^{q-1}\Delta t-|c| |X_k|^{\gamma-1}|\Delta B_k|-1-|a|\Delta t\Big|.\endaligned\]

Take $\Delta t$ small enough such that $|a|\Delta t\le 1.$ If
$|X_k|\ge\frac{2^{k+\frac{3}{q-1}}}{(|b|\Delta t)^\frac{1}{q-1}}$
and $|\Delta B_k|\le\frac{1}{2|c|}2^{(k+\frac{3}{q-1})(q-\gamma)}$, then
\[\aligned |X_{k+1}|&\ge\frac{2^{k+\frac{3}{q-1}}}{(|b|\Delta t)^\frac{1}{q-1}}
(2^{k(q-1)+3}(1-\frac{1}{2})-2)\\&
=\frac{2^{k+1+\frac{3}{q-1}}}{(|b|\Delta t)^\frac{1}{q-1}}.\endaligned\]

Thus, given that $|X_1|\ge\frac{2^{\frac{q+2}{q-1}}}{(|b|\Delta t)^\frac{1}{q-1}}$,
the event that $\{|X_k|\ge\frac{2^{k+\frac{3}{q-1}}}{(|b|\Delta t)^\frac{1}{q-1}},
\forall 1\le k\le K\}$ contains the event that $\{|\Delta B_k|\le
\frac{1}{2|c|}2^{(k+\frac{3}{q-1})(q-\gamma)},\forall 1\le k\le K\}$. So
\[P(|X_k|\ge\frac{2^{k+\frac{3}{q-1}}}{(|b|\Delta t)^\frac{1}{q-1}},\forall 1\le k\le K)
\ge\prod_{k=1}^KP(|\Delta B_k|\le\frac{1}{2|c|}2^{(k+\frac{3}{q-1})(q-\gamma)},\forall 1\le k\le K).\]

We have used the fact that $\{\Delta B_k\}$ are independent in the above inequality. But
\[\aligned P(|\Delta B_k|\ge\frac{1}{2|c|}2^{(k+\frac{3}{q-1})(q-\gamma)})&=
P(\frac{|\Delta B_k|}{\sqrt{\Delta t}}\ge\frac{ 2^{(k+\frac{3}{q-1})(q-\gamma)}}{2|c|\sqrt{\Delta t}})\\&
=\frac{2}{\sqrt{2\pi}}\int_{\frac{ 2^{(k+\frac{3}{q-1})(q-\gamma)}}{2|c|
\sqrt{\Delta t}}}^\infty e^{-\frac{x^2}{2}}dx.\endaligned\]

We can take $\Delta t$ small enough such that ${\frac{ 2^{(k+\frac{3}{q-1})
(q-\gamma)}}{2|c|\sqrt{\Delta t}}}\ge 2,$ so $x\le\frac{x^2}{2}$ for $x\ge
{\frac{ 2^{(k+\frac{3}{q-1})(q-\gamma)}}{2|c|\sqrt{\Delta t}}}$ and therefore,
\[\aligned P(|\Delta B_k|\ge\frac{1}{2|c|}2^{(k+\frac{3}{q-1})(q-\gamma)})&
\le\frac{2}{\sqrt{2\pi}}\int_{\frac{ 2^{(k+\frac{3}{q-1})(q-\gamma)}}{2|c|
\sqrt{\Delta t}}}^\infty e^{-x}dx\\&=\frac{2}{\sqrt{2\pi}}\exp\{
-{\frac{ 2^{(k+\frac{3}{q-1})(q-\gamma)}}{2|c|\sqrt{\Delta t}}}\}.\endaligned\]

Since
\[\log (1-u)\ge -2u,\quad 0<u<\frac{1}{2},\]
we have
\[\aligned\log P(|X_k|\ge\frac{2^{k+\frac{3}{q-1}}}{(|b|\Delta t)^\frac{1}{q-1}},
\forall 1\le k\le K)&\ge\sum_{k=1}^K\log (1-\exp(-{\frac{ 2^{(k+\frac{3}{q-1})(q-\gamma)}}{2|c|\sqrt{\Delta t}}}))\\&
\ge-2\sum_{k=1}^K\exp(-{\frac{ 2^{(k+\frac{3}{q-1})(q-\gamma)}}{2|c|\sqrt{\Delta t}}}).\endaligned\]

Next, by using the fact that $r^x\ge r(1\wedge\log r)x$ for any $ x\ge 1, r>1,$ we have
\[\aligned\sum_{k=1}^K\exp(-{\frac{ 2^{(k+\frac{3}{q-1})(q-\gamma)}}{2|c|
\sqrt{\Delta t}}})&=\sum_{k=1}^K\exp(-{\frac{ 2^{\frac{3(q-\gamma)}{q-1}}}{2|c|\sqrt{\Delta t}}}(2^{q-\gamma})^k)\\&
\le\sum_{k=1}^K\exp(-{\frac{ 2^{\frac{3(q-\gamma)}{q-1}}}{2|c|\sqrt{\Delta t}}}2^{q-\gamma}(1\wedge\log 2^{q-\gamma})k)\\&
\le\frac{e^{-\frac{\alpha}{\sqrt{\Delta t}}}}{1-e^{-\frac{\alpha}{\sqrt{\Delta t}}}}
\le 2e^{-\frac{\alpha}{\sqrt{\Delta t}}}\endaligned\]
for $\Delta t$ small enough, where $\alpha:=\frac{1}{2|c|}2^{\frac{(q+2)(q-\gamma)}{q-1}}(1\wedge\log 2^{q-\gamma}).$

Hence
\[\aligned\log P(|X_k|\ge\frac{2^{k+\frac{3}{q-1}}}{(|b|\Delta t)^\frac{1}{q-1}},
\forall 1\le k\le K)\ge -4e^{-\frac{\alpha}{\sqrt{\Delta t}}}.\endaligned\]

We complete the proof. $\square$

When $0<q<1, \frac{1}{2}\le\gamma<1, |b|<a,$ we also have the divergence result of the EM approximation:

\begin{Lemma}\label{divergence2}
For any $\Delta t>0$ small enough, if $|X_1|\ge r:=1+\frac{a-|b|}{2}\Delta t,$ then
there exist $k_0\ge 1$ $($depending on $\Delta t)$, A and $\alpha>0$ such that
\[\log P(|X_k|\ge r^k,\forall k\ge 1)\ge A-\frac{2e^{-k_0\alpha}}{1-e^{-k_0\alpha}}>-\infty,\]
where $A$ is finite, $\alpha:=\frac{(a-|b|)\Delta t}{2|c|}r^{1-\gamma}(1\wedge\log r^{1-\gamma}).$
\end{Lemma}

\textbf{Proof}: First, we show that
\[|X_k|\ge r^k\ \textrm{and}\ |\Delta B_k|\le\frac{(r-1)r^{k(1-\gamma)}}{|c|}\Rightarrow |X_{k+1}|\ge r^{k+1}.\]

Now
\[\aligned|X_{k+1}|&=|X_k|\Big| b|X_k|^{q-1}\Delta t+c\cdot\textrm{sgn}(X_k) |X_k|^{\gamma-1}\Delta B_k+1+a\Delta t\Big|\\&
\ge |X_k|\Big|1+a\Delta t-|b||X_k|^{q-1}\Delta t-|c| |X_k|^{\gamma-1}|\Delta B_k|\Big|\\&
\ge r^k(1+a\Delta t-|b|\Delta t-|c| r^{k(\gamma-1)}\frac{(r-1)r^{k(1-\gamma)}}{|c|})\\&
=r^k(1+2(r-1)-(r-1))=r^{k+1}.\endaligned\]

Thus, given that $|X_1|\ge r$,
the event that $\{|X_k|\ge r^k,\forall k\ge 1\}$ contains the event that
$\{|\Delta B_k|\le \frac{(r-1)r^{k(1-\gamma)}}{|c|},\forall k\ge 1\}$.

If $\frac{(r-1)r^{k(1-\gamma)}}{|c|\sqrt{\Delta t}}\ge 2$, then
\[\aligned P(|\Delta B_k|\ge\frac{(r-1)r^{k(1-\gamma)}}{|c|})&
=\frac{2}{\sqrt{2\pi}}\int_{\frac{(r-1)r^{k(1-\gamma)}}{|c|\sqrt{\Delta t}}}^\infty
e^{-\frac{x^2}{2}}dx\\&\le\frac{2}{\sqrt{2\pi}}\int_{\frac{(r-1)r^{k(1-\gamma)}}{|c|\sqrt{\Delta t}}}^\infty e^{-x}dx\\&
=\frac{2}{\sqrt{2\pi}}\exp(-\frac{(r-1)r^{k(1-\gamma)}}{|c|\sqrt{\Delta t}}).\endaligned\]

We can choose $k_0$ be the smallest $k$ such that $\frac{(r-1)r^{k(1-\gamma)}}{|c|\sqrt{\Delta t}}\ge 2$
(note that since $r>1$, such $k_0$ always exists).

On the other hand,
\[\aligned &\sum_{k=k_0}^\infty\log(1-\frac{2}{\sqrt{2\pi}}\exp(-\frac{(r-1)r^{k(1-\gamma)}}{|c|\sqrt{\Delta t}}))\\&
\ge -2\sum_{k=k_0}^\infty\exp(-\frac{(r-1)r^{k(1-\gamma)}}{|c|\sqrt{\Delta t}})\\&
\ge -2\sum_{k=k_0}^\infty\exp(-k\times\frac{(r-1)r^{1-\gamma}(1\wedge\log r^{1-\gamma})}{|c|\sqrt{\Delta t}})\\&
=-\frac{2e^{-k_0\alpha}}{1-e^{-k_0\alpha}}>-\infty.\endaligned\]

So $\prod_{k=1}^\infty P(|\Delta B_k|\le\frac{(r-1)r^{k(1-\gamma)}}{|c|})$ is well defined and therefore
\[P(|X_k|\ge r^k, \forall\ k\ge 1)\ge\prod_{k=1}^\infty P(|\Delta B_k|\le\frac{(r-1)r^{k(1-\gamma)}}{|c|}).\]

Then as in proof of Lemma \ref{divergence1}, we have
\[\aligned \log P(|X_k|\ge r^k, \forall\ k\ge 1)&\ge
A-\frac{2e^{-k_0\alpha}}{1-e^{-k_0\alpha}}>-\infty,\endaligned\]
where
\[A=\sum_{k=1}^{k_0-1} \log P(|\Delta B_k|\le\frac{(r-1)r^{k(1-\gamma)}}{|c|}),\]
\[\alpha=\frac{(r-1)r^{1-\gamma}(1\wedge\log r^{1-\gamma})}{|c|\sqrt{\Delta t}}.\]

We complete the proof. $\square$

Let us give an example to show that the $\theta$-EM scheme
($\frac{1}{2}<\theta\le 1$) is exponentially stable while EM scheme
is not.

\textbf{Example 1:}

Consider the following one dimensional stochastic differential equation,
\begin{equation}\label{sde2}
dX_t=(aX_t+b|X_t|^{2\gamma-2} X_t)dt+c|X_t|^\gamma dB_t,
\end{equation}
where $\gamma>1$ $a<0$ and $2b+c^2\le 0$.

It is clear that both of the coefficients are locally Lipschitz
continuous. Thus SDE (\ref{sde2}) has a unique global solution.

By Lemma \ref{divergence1}, since $2\gamma-1>\gamma>1$, we know that
when we choose the step size $\Delta t$ small enough and the initial
value $X_1$ suitably, the classical EM scheme is divergent with a
positive probability. Now let us consider the exponential stability
of $\theta$-EM scheme.

The corresponding $\theta$-EM scheme of (\ref{sde2}) is
\begin{equation}\label{STM2}
X_{k+1}=X_k+[(1-\theta)X_k(a+b|X_k|^{2\gamma-2})+\theta X_{k+1}
(a+b|X_{k+1}|^{2\gamma-2})]\Delta t+c|X_{k}|^\gamma \Delta B_k,
\end{equation}

Notice that in our case $g(x)=c|x|^\gamma$ does not satisfy the
linear growth condition. Therefore, the stability results in
\cite{HMS,HMY,PDM,WMS} as well as \cite{CW} for the moment as well
as almost sure exponential stability of the backward EM scheme case
($\theta=1$) can not be used here.

On the other hand, since in this case $f(x)=ax+b|x|^{2\gamma-2} x,
g(x)=c|x|^\gamma,$ it is obvious that
$$2\langle x,f(x)\rangle+|g(x)|^2=2ax^2+(2b+c^2)|x|^{2\gamma}\le 2ax^2.$$

Since $a<0$, then condition (\ref{c2}) holds for $C=-2a$. Moreover,
$$\langle x-y,f(x)-f(y)\rangle=a(x-y)^2+b(x-y)(|x|^{2\gamma-2}x-|y|^{2\gamma-2}y).$$

Since
$$(x-y)(|x|^{2\gamma-2}x-|y|^{2\gamma-2}y)\ge0$$
holds for $\forall x,y\in\mathbb{R},$ it follows that
$$\langle x-y,f(x)-f(y)\rangle\le a(x-y)^2.$$

We have used the fact that $b<0$ here. Thus conditions (\ref{c2})
and (\ref{c3}) hold. By Theorem \ref{exponential}, we know that, for
any $0<\varepsilon<1$, the $\theta$-EM ($\frac{1}{2}<\theta\le 1$)
scheme (\ref{STM2}) of the corresponding SDE (\ref{sde2}) is mean
square exponentially stable with with Lyapunov exponent no greater
than $2a(1-\varepsilon)$ and almost surely exponentially stable with
Lyapunov exponent no greater than $a(1-\varepsilon)$ if $\Delta t$
is small enough.

For the polynomial stability, we consider the following example.

\textbf{Example 2:}

Now let us consider the following scalar stochastic differential equation,
\begin{equation}\label{sde3}
dX_t=\frac{-(1+t)^{\frac{1}{2}}|X_t|^{2\gamma-2}X_t-2K_1X_t}{2(1+t)}dt+
\sqrt{\frac{|X_t|^{2\gamma}}{(1+t)^\frac{1}{2}}+\frac{C}{(1+t)^{K_1}}}
dB_t,
\end{equation}
where $C>0$, $K_1>1$, $\gamma\ge1$ are constants.

Since in this case
$$f(x,t)=\frac{-(1+t)^\frac{1}{2}|x|^{2\gamma-2}x-2K_1x}{2(1+t)},\quad g(x,t)=
\sqrt{\frac{|x|^{2\gamma}}{(1+t)^\frac{1}{2}}+\frac{C}{(1+t)^{K_1}}},$$
It is clear that both of the coefficients are locally Lipschitz
continuous. Moreover, it is easy to verify that
$$2\langle x,f(x,t)\rangle+|g(x,t)|^2\le C(1+t)^{-K_1}- K_1(1+t)^{-1}|x|^2,$$
and
$$\langle x-y,f(x,t)-f(y,t)\rangle\le 0\le L|x-y|^2.$$

Thus conditions (\ref{c1}) and (\ref{c3}) hold. Therefore, SDE
(\ref{sde3}) has a unique global solution. If $\gamma>1,$ then by
Theorem \ref{polynomial}, for any $0<\varepsilon<K_1-1,$ the
$\theta$-EM ($\frac{1}{2}<\theta\le 1$) solution of (\ref{sde3})
satisfies the polynomial stability (with rate no great than
$-(K_1-1-\varepsilon)$) for $\Delta t$ small enough. If $\gamma=1,$
it is obvious that $f$ also satisfies the linear growth condition
(\ref{growth}) (condition (2.4) in \cite{LFM} failed in this case),
then by Theorem \ref{polynomial1}, the $\theta$-EM
($0\le\theta\le\frac{1}{2}$) solution of (\ref{sde3}) satisfies the
polynomial stability for $\Delta t$ small enough. However, since the
coefficient $g(x,t)$ is not bounded with respect to $x$, we can not
apply Theorem 3.1 and Theorem 3.5 in \cite{LFM} to get the
polynomial stability of the classical EM scheme and back EM scheme
respectively.

\textbf{Acknowledgement:} The second named author would like to
thank Professor Chenggui Yuan for useful discussions and suggestions
during his visit to Swansea University.

\end{document}